\documentclass{article}
\usepackage{amsmath,amsthm,amssymb,psfig}

\newtheorem{lemma}{\bf Lemma}
\newtheorem{theorem}{\bf Theorem}

\newtheorem{remark}{\bf Remark}

\begin{document}
\title{Uncertainty estimates and $L_2$ bounds for the 
Kuramoto-Sivashinsky equation}
\author{Jared C. Bronski\footnote{Department of Mathematics,University of Illinois Urbana-Champaign, 1409 W. Green St,
Urbana IL 61801}\\ Tom Gambill${}^*$}

\maketitle

\noindent
\begin{abstract} We consider the Kuramoto-Sivashinsky (KS) equation in one 
dimension with periodic boundary conditions. We apply a Lyapunov function 
argument similar to the one first introduced by Nicolaenko, Scheurer, and 
Temam \cite{NICO}, and later improved by Collet, Eckmann, Epstein and 
Stubbe\cite{CEES}, and Goodman \cite{Good} to prove that 
$\limsup_{t\rightarrow \infty} |\!|u|\!|_2 \le C L^{\frac{3}{2}}.$ 
This result is slightly
weaker than that recently announced by Giacomelli and Otto \cite{GO}, but 
applies in the presence of an additional linear destabilizing term. We further 
show that for a large class of functions
$\phi_x$ the exponent $\frac{3}{2}$ is the best possible from this line of argument. 
Further, this result together  
with a result of  Molinet\cite{Mol} gives an improved estimate for $L_2$ 
boundedness of the Kuramoto-Sivashinsky equation in thin rectangular 
domains in two spatial dimensions.

\end{abstract}

\section{INTRODUCTION}
\subsection{Background}
The Kuramoto-Sivashinsky (KS) equation
\[
u_t = -u_{xxxx} - u_{xx} + u u_x ~~~~~ \int u(x,0)dx = 0 
\]
arises as a model of certain
hydrodynamic problems, most notably the propagation of flame fronts\cite{Siv}. 
The KS equation is interesting
mathematically because the linearization about the zero state
\[
u_t = -u_{xxxx} - u_{xx}
\]
subject to periodic boundary conditions on $[-L,L]$ has a large
number ($O(L/\pi)$) of exponentially growing modes. The growth of
these modes corresponds, in the combustion problem,
to the development of nontrivial structures. In addition to its 
importance as a model for flame fronts\cite{Siv} and 
phase turbulence\cite{KS} and plasmas\cite{LMRT} the 
KS equation has become 
one of the canonical models for 
spatio-temporal chaos in $1\!+\!1$ dimensions\cite{JOHNSON,KEVREKIDIS,Man}. 

Nicolaenko, Scheurer and Temam\cite{NICO} gave the first long-time 
boundedness result for the Kuramoto-Sivashinsky equation, showing 
 that $\limsup_{t \rightarrow \infty} |\!|u|\!|_2 
\le C L^{\frac{5}{2}}$ for odd initial data, as well as showing 
that bounds on the $L_2$ norm imply bounds on the dimension of the attractor. 
The $L_2$ estimate was improved by 
Collet, Eckmann, Epstein and 
Stubbe\cite{CEES} who  extended it to any mean-zero initial data and improved the exponent from $\frac{5}{2}$ to $\frac{8}{5}$, 
and by Goodman\cite{Good}, who extended it to 
any mean-zero initial data with the same exponent. All of these papers use a 
version of the original argument of Nicolaenko, Scheurer and Temam, namely to 
establish that the function $|\!|u-\phi|\!|_2^2$ is a Lyapnuov function for 
an appropriately chosen $\phi$ and $|\!|u|\!|_2$ sufficiently large. There are 
also two bounds which do not fit into this Lyapunov function framework, 
that of Ilyashenko\cite{Ilyashenko}, and the very recent paper of Otto 
and Giacomelli\cite{GO}. The latter, which treats the KS equation as a 
perturbation of the Burger's equation, is currently the best estimate, establishing that
\[
\limsup_{t \rightarrow \infty} |\!|u|\!|_2 = o(L^{\frac{3}{2}}). 
\]   

In this paper we give an elementary argument of the Lyapunov function 
type which establishes the slightly weaker result
\[
\limsup_{t \rightarrow \infty} |\!|u|\!|_2 = O(L^{\frac{3}{2}}).
\]
Our proof applies equally to the destabilized Kuramoto-Sivashinsky (dKS)
equation:
\[
u_t = -u_{xxxx} - u_{xx} + \gamma  u + u u_x ~~~~~~\gamma > 0 ~~~~~ \int u(x,0)dx = 0. 
\]
It was shown by Wittenberg\cite{WITTENBERG} that this equation 
has stationary solutions which satisfy $|\!|u|\!| \propto L^{\frac{3}{2}}.$ 
Since a Lyapunov function argument for the KS equation also applies 
to the dKS equation (for sufficiently small $\gamma$) Wittenberg 
argued that $\frac{3}{2}$ is the best exponent that one can expect 
from the Lyapunov function approach. This paper completes this circle 
of ideas, by showing that this exponent can actually be achieved. 
 
We also give a independent scaling argument that motivates the choice of 
$\phi_x$ and shows that for a large class 
of potentials $\frac{3}{2}$ is the best exponent possible from 
this line of argument. This argument is useful because it makes clear 
the physical basis of the scaling, and is potentially applicable to 
other equations. 

\subsection{Fundamental Lemmas}
We begin by stating two basic lemmas which form the core of the Lyapunov
function argument. These lemmas are basically equivalent to equations 2.11,12 
and 3.10-3.12 in \cite{NICO}, or analogous results in
\cite{CEES},\cite{Good}. It is worth noting that similar ideas of 
considerably greater generality have been used by Constantin and Doering
to establish bounds on energy dissipation in fluids, and generally 
go by the name `background flow method.'\cite{CD1,CD2} 

\begin{lemma} \label{thelemma}
Given $u=u(x,t) \in { L_2} [-L,L]$ and $\phi = \phi(x) \in {L_2} [-L,L]$ 
satisfying the following inequality:
\begin{eqnarray}
 & \frac{\partial}{\partial t} |\!|u-\phi|\!|_2^2 \leq -\lambda|\!|u|\!|_2^2+ M^2~ \label{main}
\end{eqnarray}
for some constants $\lambda > 0$ and $M$, then $B\left(0,R^{**}\right)$ , 
the ball of radius $R^{**}$ centered about the origin, is an attracting 
region, where the radius $R^{**}$ is given by 
\begin{eqnarray}
R^{**} = \sqrt{2|\!|\phi|\!|_2^2+\frac{2M^2}{\lambda}}+|\!|\phi|\!|_2.
 \label{rstarstar} 
\end{eqnarray}
\end{lemma}
To show this, let $|\!|u|\!|_2 \geq R^{**}$ then
\begin{eqnarray*}
|\!|u-\phi|\!|_2\geq|\!|u|\!|_2-|\!|\phi|\!|_2\geq\sqrt{2|\!|\phi|\!|_2^2-\frac{2M^2}{\lambda}}\end{eqnarray*}
The parallelogram law implies 
\[-\lambda|\!|u-\phi|\!|_2^2\geq -2\lambda|\!|u|\!|_2^2-2\lambda|\!|\phi|\!|_2^2\] 
which in turn gives
\begin{eqnarray*}
 & \frac{\partial}{\partial t} |\!|u-\phi|\!|_2^2 \leq \lambda|\!|\phi|\!|_2^2+
M^2 -\frac{\lambda}{2}|\!|u-\phi|\!|_2^2. 
\end{eqnarray*}
If we apply the obvious Gronwall estimate to the above inequality 
it is apparent that $B(\phi,R^*)$, the ball of radius 
$R^*$ centered about $\phi$ is exponentially attracting, with ${R^*}^2 = |\!|\phi|\!|^2 + \frac{2 M^2}{\lambda}$. The triangle inequality implies 
$B(\phi,R^*) \subset B(0,R^{**}).$\\
$\clubsuit$

\begin{lemma}
For any $\phi\in \dot H^2_{\rm per}$ and $u(x,t)$ solving the 
Kuramoto Sivashinsky equation we have (after some rescaling) the inequality
\[
\frac{1}{128}\frac{\partial}{\partial t}\int_{-64L}^{64L}(u - 16 \phi)^2  \leq
8\bigl(\int_{-64 L}^{64 L}  u_{y}^2- u_{yy}^2-\tilde \phi_y u^2\bigr)+\int_{-64 L}^{64 L} 16\tilde \phi_y^2 +2(16^3)\tilde \phi_{yy}^2~~~~~
\]
\end{lemma} 
A straightforward calculation gives
$$
 \frac{1}{2}\frac{\partial}{\partial t}|\!|u - \phi|\!|_2^2 = \int_{-L}^{L}u_t(u-\phi) = \int (-u_{xx}-u_{xxxx}-uu_x)(u-\phi). \label{lyaponov}  
$$
After integrating by parts and applying periodic boundary conditions this becomes
\begin{eqnarray*}
 \frac{1}{2}\frac{\partial}{\partial t}|\!|u - \phi|\!|_2^2 &=&
 \int u_{x}^2-u_{xx}^2-\phi_xu_x+\phi_{xx}u_{xx}-\frac{1}{2}\phi_xu^2.
\end{eqnarray*}
Applying the Cauchy-Schwartz inequality in the form 
$<\!\!f,g\!\!> \leq \frac{p}{2}\!<\!\!f,f\!\!> + \frac{1}{2p}\!<\!\!g,g\!\!>$ gives
\begin{eqnarray*}
 \frac{1}{2}\frac{\partial}{\partial t}|\!|u - \phi|\!|_2^2 &\leq&
\int (1+\frac{1}{2p})u_{x}^2+(\frac{1}{2q}-1)u_{xx}^2+\frac{p}{2}\phi_x^2 +\frac{q}{2}\phi_{xx}^2-\frac{1}{2}\phi_xu^2.
\end{eqnarray*}
If we then make the substitution $\phi = \gamma \tilde \phi, ~y =\beta x,$
we find that 
\begin{eqnarray*}
\frac{1}{128}\frac{\partial}{\partial t}\int_{-\beta L}^{\beta L}(u - 16 \phi)^2 dy
\le \int_{-\beta L}^{\beta L} \frac{1+2p}{2p}\beta  u_{y}^2+
\frac{1-2q}{2q}\beta ^3 u_{yy}^2+\frac{p}{2}\beta  \gamma ^2\tilde \phi_y^2 
+\frac{q}{2}\beta ^3 \gamma ^2 \tilde \phi_{yy}^2-\frac{1}{2} \gamma  
\tilde \phi_y u^2 dy~.
\end{eqnarray*}
Finally, taking $p = \frac{1}{2}, q = 1, \beta = 64, \gamma = 16 $ we get 
\[
\frac{1}{128}\frac{\partial}{\partial t}\int_{-64L}^{64L}(u - 16\phi)^2 dy 
\leq 8\int_{-64 L}^{64 L} \bigl(u_{y}^2- u_{yy}^2-\tilde \phi_y u^2\bigr) dy 
+\int_{-64 L}^{64 L} 16\tilde \phi_y^2 +2(16^3)\tilde \phi_{yy}^2,
\]
as claimed.\\ $\clubsuit$

\begin{remark}
Since we are only concerned with the scaling of the 
estimates with $L$, and not with the actual constants, we will henceforth 
drop the tilde and replace $\beta L$ with $L$ and $y$ with $x$. 
The preceeding 
lemmas show that, if we can construct $\phi\in \dot H^2_{per}$ such that 
the coercivity estimate 
\[
<\!\!u,Ku\!\!> = \int u_{xx}^2 - u_x^2 + \phi_x u^2 dx > \delta|\!|u|\!|_2^2 > 0
\]
for some $\delta$ {\it independent of $L$}, then we get an estimate of the form 
\[
\limsup_{t\rightarrow\infty} |\!|u|\!|_2 \leq  R^{**} = \sqrt{c_1 |\!|\phi|\!|_2^2 + c_2 |\!|\phi_x|\!|_2^2 + c_3 |\!|\phi_{xx}|\!|_2^2 } + c_4 |\!|\phi|\!|_2^2 .
\]
Since it is clear that that $R^{**}$ is comparable to the 
$H^2$ norm, $c |\!|\phi|\!|_{H^2} \le R^{**} \le C |\!|\phi|\!|_{H^2}$, 
we will write this in the form
\[
\limsup_{t\rightarrow\infty} |\!|u|\!|_2 \leq c |\!|\phi|\!|_{H^2}.
\]   
\end{remark}

Our argument proceeds along the same lines as that of
 previous papers, but with 
a better construction of the function $\phi_x$, which leads to better 
exponents than in previous Lyapunov function arguments. 
We first construct a $\phi_x$ such that the operator $K$ defined above 
is positive definite for  $u$ satisfying Dirichlet boundary
condition. This establishes an $L_2$ bound for odd solutions of the KS
equation, which are preserved under the flow. This result can be extended to 
all mean-zero data by allowing a time-dependent $\phi_x$ which translates 
under a kind of gradient-flow dynamics, as was first done by Collet, Eckmann, 
Epstein and Stubbe. Our construction is done in real space, unlike 
the constructions of Nicolaenko, Scheurer and Temam, and Collet, Eckmann, 
Epstein and Stubbe, where the construction uses a clever cancellation 
in Fourier space. We feel that the real space construction clarifies the 
role of uncertainty estimates in determining the positivity of the 
operator, and makes it easier to see what the optimal scaling should be.

\section{MAIN RESULTS}
\subsection{Scaling, Uncertainly, and Bounds}

In this section we present our main results. 
We begin with a discussion of the role of scaling and uncertainty
in determining exponents in the Lyapunov function approach to 
proving boundedness of the KS equation. The construction of a suitable 
function $\phi_x$ can be viewed as a competition between kinetic energy 
and potential energy terms in the operator, and relatively simple 
scaling arguments make it clear how the function should scale 
with $L$, the length of the interval. With this as motivation, we 
proceed to prove that a suitable Lyapunov function with the critical 
scaling exponents can be constructed. The main technical tools will be 
a Hardy-type inequality, which allows us to derive a lower bound on a
second order kinetic energy term with a Dirichlet boundary condition in 
terms of a standard first-order kinetic energy term, together with a 
lower bound on a Schrodinger operator in terms of a finite dimensional 
quadratic form.   

As outlined in the previous section the basic strategy is to choose a 
periodic function $\phi_x$ of zero mean such that the following 
quadratic form is coercive,
\[
<\!\!u, K u\!\!>= \int u_{xx}^2 - u_x^2 + \phi_x u^2 \ge \delta |\!|u|\!|^2, 
\]
for $u$ satisfying Dirichlet boundary conditions and some positive 
$\delta$ independent of $L$.  Doing so gives a bound on the 
radius of the attracting ball in $ L_2$ (for odd solutions)
which scales like $|\!|\phi|\!|_{H^2}$. A sketch of the 
strategy which has previously 
been followed for constructing such a potential\cite{NICO,CEES} is this: One 
constructs a potential in Fourier space which is constant for a 
large range of wavenumbers, and decaying thereafter. For such a 
potential the quadratic form above looks like a diagonal piece 
plus a piece which is supported only at very large wavenumbers.
One then uses the rapid growth of the dispersion relation $k^4 - k^2$ 
to show diagonal dominance and thus positivity of the operator. 
The real space translation of this strategy is as follows:  
One constructs a potential that is large and negative on a small set near 
the origin, and positive on the rest of the interval, in such a way that 
the potential is mean zero. 
The uncertainty principle, together with the Dirichlet boundary condition,  
implies that little of the mass of the ground state 
can be concentrated in the small region where the potential is negative, 
so that the net effect is to produce a positive ground state 
eigenvalue. In this section we present a scaling argument which 
suggests the best possible estimate one can get of this form. In the 
next section we show that this estimate can actually be achieved. 

We will assume for simplicity of discussion that for $x \in (-L,L)$ 
the function $\phi_x$ takes the following form, 
$$
\phi_x = \gamma L^{c_2 - c_1 - 1} +  L^{c_2} \tilde q(x L^{c_1})~~~~~~~~ \gamma,c_{1,2}>0 \label{phiform}
$$
where $\gamma$ is a constant and $\tilde q$ is compactly supported $C^2$ 
function, with $\phi_x$ extended to a $2L$ 
periodic function in the usual way. The functions constructed by 
Goodman and in the current paper are exactly of this form, while
the functions $\phi_x$ constructed by other authors 
can all be written as a sum of a constant and a function which is rapidly 
decaying
away from the origin, though not necessarily compactly supported. The present 
discussion can be extended to such potentials under some 
very mild technical assumptions. For details see the Ph.D. thesis of one of 
the authors\cite{Tom}. For reference the functions constructed by Temam 
et. al and by Goodman have scaling exponents $c_1=1,c_2=2$ 
while the function constructed by Collet et. al. has scaling exponents  
$c_1=\frac{2}{5}, c_2 = \frac{7}{5}$.

Our first observation is that the operator cannot be positive for 
$c_2-c_1-1<0$. This follows from a straightforward test-function argument 
using a delocalized 
test function such as $u = L^{-\frac{1}{2}} \sin(\frac{k \pi x}{L})$,
for suitably chosen $k$.  
Next, if one makes the rescaling $y = L^{c_1} x,$ the quadratic form becomes
\[
L^{3 c_1}\left(\int_{-L^{1+c_1}}^{L^{1+c_1}} u_{yy}^2 - L^{-2 c_1} u_y^2 + L^{c_2-4c_1}\tilde q(y) u^2(y) dy \right) + \gamma L^{c_2-c_1-1} |\!|u|\!|_2^2.  
\]

Note the prefactor of $L^{c_2 - 4 c_1}$ in front of the potential $\tilde q$.
Motivated by this, we refer to potentials for which $c_2 - 4 c_1<0$ as 
weak potentials, those for which  $c_2 - 4c_1> 0$ as strong potentials, 
and those for which $c_2 = 4 c_1$ as critical potentials. Strong potentials 
are those for which the potential energy term dominates the low-lying  
eigenvalues in the limit $L \rightarrow \infty$, while weak potentials are 
those for which the kinetic energy 
term dominates. All of the potentials constructed in previous papers are 
weak potentials, with the potential in \cite{CEES} being closest to critical.

It is clear from another simple test-function 
argument, this time with a localized test-function, 
 that in the case of a strong potential the operator 
$K$ again cannot be positive.
Simply taking a compactly supported test function whose support 
is contained in a region where $\tilde q<0$ gives an estimate of the 
following form 
$$
\lambda_0(K) \le - C L^{c_2 - c_1} + O(L^{3c_1}, L^{c_1}, L^{c_2-c_1-1}). \label{testfn2}
$$
where the three error terms come from the $u_{xx}^2, u_x^2$ and $u^2$ terms respectively. 
A simple calculation shows that the $H^2$ norm of the potential 
$\phi_x$ is bounded below by  
$$
|\!|\phi|\!|_{H^2} \ge |\!|\phi_{xx}|\!|_2 = O(L^{c_2+\frac{c_1}{2}}).
$$
Thus the best estimate possible for a $\phi_x$ of this form 
is given by the solution to the constrained minimization problem
\begin{eqnarray}
{\rm minimize} ~~ c_2 &+& \frac{c_1}{2} ~~ {\rm subject~to} \label{exponent} \\
c_2 &\ge& c_1 + 1 \label{constraints1}\\
c_2 &\le& 4 c_1 \label{constraints2}
\end{eqnarray}

It is easy to check that the solution to this constrained minimization 
problem is given by 
$$
c_2 = \frac{4}{3}~~~~ c_1 = \frac{1}{3}.
$$ 
A cartoon of such a potential is shown in figure(1). 
A potential of this form would give an estimate of the radius of the attractor 
that scales like $L^{\frac{4}{3} + \frac{1}{2}\cdot \frac{1}{3}} = L^\frac{3}{2}.$ We show in this paper that these critical exponents 
can actually be achieved. 

It is worth noting that the scaling given above for the 
potential function $\phi_x$ is exactly the same as that for the 
viscous shock profile for the destabilized KS equation
constructed by Wittenberg\cite{WITTENBERG}. This is not too 
surprising, since the same Lyapunov function argument 
applies to the destabilized KS equation, and one might 
reasonably expect 
that the Lyapunov function should look something like the 
steady state. In fact, it is easy to check that if $\phi_\gamma$ is a 
stationary solution to the destabilized KS equation, and the 
linearized operator is negative semi-definite, then $|\!|u -\phi|\!|^2$ 
is a Lyapunov function for sufficiently large $|\!|u|\!|.$

\begin{figure}
\begin{center}
\begin{tabular}{c}
\psfig{file=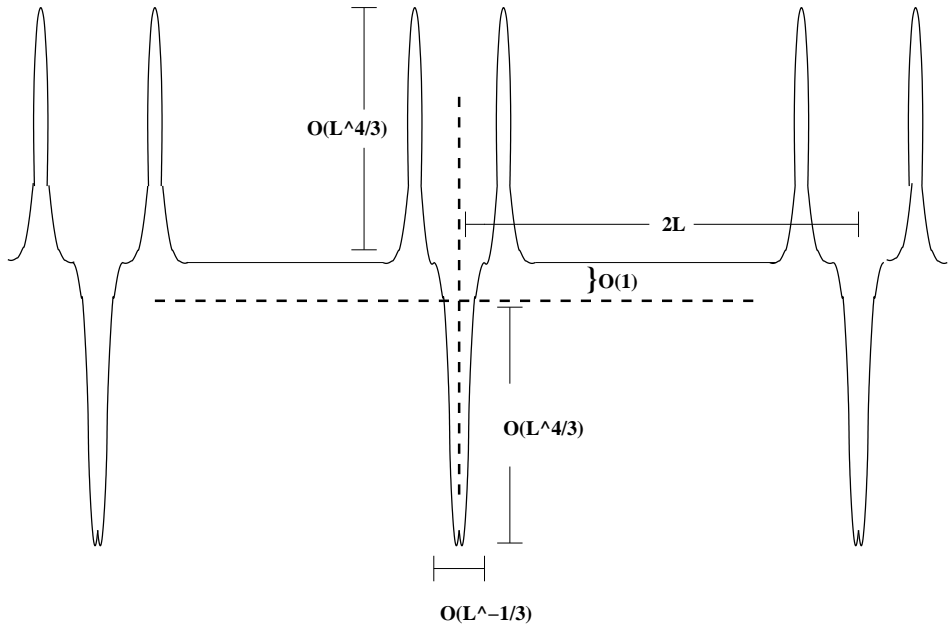,height=3.25in}
\end{tabular}
\end{center}
\caption{A cartoon of the function $\phi_x$ constructed in this paper. }
\label{fig:pot}
\end{figure}

This calculation makes clear the construction of a Lyapunov function 
for this problem involves a competition between the kinetic energy terms
in the functional (which scale with the width of the potential) and the 
potential energy terms, which scale with the height of the potential. 
In particular it should be clear from this calculation that the scaling 
depends crucially 
on the order of the operator. In particular one expects {\em different} 
critical exponents for a second order operator, since the kinetic energy 
term is less effective at small scales than the analogous term 
for a fourth order operator.  
Along these lines we make a couple of other comments, directed 
specifically at the papers of Nicolaenko, Scheurer and Temam and
of Goodman.  In \cite{NICO} and \cite{Good}
they make an additional simplification by using the inequality
$\partial_{xxxx} + \partial_{xx} \ge -\partial_{xx} + 1$ to bound the 
operator $K$ from below by a standard second order Schrodinger 
operator $\tilde K$:
\[
<\!u,\tilde K u\!> = \int u_{x}^2 - u^2 + \phi_x u^2 ~\le~ <\!u,K u\!> = \int u_{xx}^2 - u_x^2 + \phi_x u^2
\] 
If one carries out the same scaling analysis presented above 
for this quadratic form one finds that the exponents $c_1,c_2$ defined above must satisfy the inequalities
\begin{eqnarray}
c_2 &\ge& c_1 + 1 \\
c_2 &\le& 2 c_1 ,
\end{eqnarray}
along with the same estimate of the penalty term due to the potential:
\[
|\!|\phi|\!|_{H^2} \ge |\!| \phi_{xx} |\!|_2 = O(L^{c_2 + \frac{c_1}{2}}).
\]
Carrying out this minimization problem gives the critical exponents
$c_1 = 1,c_2=2$, giving an estimate of the radius of the attracting 
ball of $R = C L^\frac{5}{2}.$ Thus the estimates in \cite{NICO} and 
\cite{Good}
are the best possible for potentials of this form and the 
second order operator $\tilde K$.

This calculation is, in a sense, complementary to Lieb-Thirring type 
inequalities.
In Lieb-Thirring inequalities one attempts to maximize some measure of the 
negative part of the spectrum of an operator over all potentials with 
a fixed norm. In the Lyapunov function method one would like to minimize
the negative part of the spectrum, to obtain a positive operator. 
Unfortunately, there appear to be only few results of this sort for operators 
of higher order than second (see, however, the work of Tadjbakhsh and Keller\cite{TK} 
and Cox and Overton\cite{CO} on the optimal shape of columns).

As a sidenote we remark that the Burgers-Sivashinsky equation, being 
second order, provides an amusing application of this point of view. 
An analysis of the Burgers-Sivashinsky (BS) equation 
$$
u_t = u_{xx} + u + u u_x 
$$
by the Lyapunov function method 
leads to a bound of the form 
\[
\frac{d}{dt} |\!|u-\phi|\!|^2 \le -\lambda_0 |\!|u|\!|^2 + c |\!|\phi_x|\!|^2_2 
\]
where $\lambda_0$ is the smallest eigenvalue of the operator
\[
H u = -u_{xx} + \phi_x u
\]
One can choose the potential $\phi_x$ in an `optimal' way by maximizing 
the ground state eigenvalue $\lambda_0(\phi_x)$ as a functional of 
$\phi_x$ subject to the constraints $\int \phi_x^2 dx =
{\rm constant}.$ Formally one has 
\begin{eqnarray}
\max_{\phi_x} \lambda_0(\phi_x)|\!|u|\!|_2^2 - \mu \int\phi_x^2 &=& \\
\max_{\phi_x} \min_u \int u_x^2 + \phi_x u^2 - \mu \phi_x^2  &=& \\
\min_u \max_{\phi_x} \int u_x^2 + \phi_x u^2 - \mu \phi_x^2  
\end{eqnarray}
where $\mu$ is the Lagrange multipliers associated with the 
constraint. Of course the replacement of the $\max\min$ by $\min\max$
needs justification, but this can be done fairly easily in the second order case.
 Solving the maximization over $\phi_x$ leads to 
$\phi_x = \frac{1}{2 \mu} \left( u^2 - |\!|u|\!|_2^2/(2L)\right)$ and thus to 
the functional
\[ 
\min_u  \int u_x^2 + \frac{1}{4 \mu}(u^2 - |\!|u|\!|_2^2/(2L))^2 
\]
The solution of this minimization problem can be expressed in closed form 
in terms of Jacobi elliptic functions\cite{Tom}.

\subsection{Proof of Main Results}
In this section we show that the critical exponents $c_1 = \frac{1}{3},c_2=\frac{4}{3}$ given by the solution to 
 Eqs (\ref{exponent},\ref{constraints1},\ref{constraints2}) 
can actually be achieved. Our main techniques are a higher-order 
analog of the Hardy inequality together with an elementary 
uncertainty estimate, which allow us to bound the quadratic form from 
below by a finite-dimensional one. Our main result can be stated as follows

\indent\begin{theorem}

$\exists$ a $2L$ periodic potential function $\phi_x$ such that 
\[
\int u_{xx}^2 - u_x^2 + \phi_x u^2 \ge \frac{1}{4} \int u_{xx}^2 + u^2 
\]
$\forall u \in C^2[-L,L]$ with $u(0)=0$. Further we have the 
estimate $|\!| \phi|\!|_{H^2}\le c L^{\frac{3}{2}}.$

\end{theorem}

{\bf Note:} The above theorem requires only a single Dirichlet 
boundary condition at the origin, and thus applies to many different  
densely defined domains for the operator $K = \partial_{xxxx} + \partial_{xx}
+ \phi_x.$  We will primarily be interested in the domain of odd, periodic 
flows, which is preserved under the KS flow. 

The proof of this theorem is presented as a series of simple lemmas.
The first lemma is an uncertainty inequality that allows us to estimate 
the second order  kinetic energy term by a first order kinetic energy 
which is more analytically tractable. This lemma has an advantage over 
the standard Poincare inequality, since it is in some sense `local' and 
does not scale badly with large intervals. The downside is that we `use up'
the Dirichlet boundary condition, so we are forced to consider a larger 
domain for $v$. 

\indent \begin{lemma} \label{hardylemma}
Suppose that $u \in C^\infty$ with $u(0)=0$. Then, if $v(y) = \frac{u(y)}{y}$
we have the inequality
\begin{eqnarray*}
 \int_{-a}^{a} \frac{1}{4}u_{yy}^2(y) \ge \int _{-a}^a 
\frac{1}{2} \left(\frac{u(y)}{y}\right)_y^2 = 
\int_{-a}^{a} \frac{1}{2}v_y^2 dy  
\end{eqnarray*}
\end{lemma}
\textbf{Proof:} 
Since $u(y) = y v(y)$ we have 
$u_{yy} = yv_{yy}+2v_y$ and after substitution into the above integral we
get
\begin{eqnarray*}
\int_{-a}^{a} u^2_{yy} dy & = & \int_{-a}^{a} (yv_{yy} + 2v_y)^2 dy \\
& = & \int_{-a}^{a} y^2v_{yy}^2 dy + 2\int_{-a}^{a} v^2_y dy + 2 a \left(v^2(a) + v^2(-a)\right)\ge 2\int_{-a}^{a} v^2_y.
\end{eqnarray*}
Here we have integrated by parts once and used the fact that the term 
$2yv^2_y\vert_{-a}^{a}= 2a(v^2(a)+ v^2(-a))$ 
is positive. \\
$\clubsuit$

\begin{remark}
 This is essentially a higher order analog of the Hardy 
inequality, 
which says that for $F\in C^1, F(0) = 0$ one has the estimate 
\[
\int |F_x|^2 dx \ge \frac{1}{4}\int \frac{|F(x)|^2}{x^2} dx  
\]
\end{remark}

In the next lemma we show the positivity of an operator which we will 
later use to bound then operator $K$ from below.

\indent \begin{lemma} \label{PClemma}
Define a piecewise constant compactly supported function  $Q(y)$ as follows:   
\[ Q(y)  = \left\{\begin{array}{cl}
                  -q_0, & \mbox{when $0 \le |y| \le \frac{a}{2}$} \\
		   q_1 &   \mbox{when $ \frac{a}{2} < |y| \le a$} \\
                   0 &   \mbox{when $  a < |y| $}
		 \end{array}\right. ,   
\]
with $a,q_0,q_1$ positive constants satisfying the inequalities
\begin{eqnarray} 
q_0 a^2 &<& 1 \\
 q_1 &>& \frac{q_0}{1 - a^2 q_0}
\end{eqnarray}
Then for all $v \in H^1$ we have 
\[\int \frac{1}{2}v^2_y + Q(y)v^2 dy > 0 ~.\]
\end{lemma}
\textbf{Proof:} 
We will show that $\int_{0}^a \frac{1}{2}v^2_y + Q(y)v^2 dy >0 $ for any 
$v \in H^1.$ Since $Q$ is even the same argument holds for 
$\int_{-a}^0 \frac{1}{2}v^2_y + Q(y)v^2 dy.$ Obviously the integral 
over $|y| > a$ is positive (since $Q$ is zero) and can thus be neglected. 
Note that we are not assuming 
any particular boundary conditions on $v$. 

For any $v\in H^1$ and any two points $y_1$ and $y_2$ 
we have the elementary uncertainty inequality
\[\int_{y_1}^{y_2}v^2_y \ge \frac{(v(y_1)-v(y_2))^2}{y_2-y_1}~.\]

Since $v\in H^1$ is continuous we can `sample' $v(y)$ at three 
location $y_0 = 0,~ y_1 \in (0,a/2), y_2\in(a/2,a)$ 
defined as follows
\begin{eqnarray*}
  v(y_0) &=&  v_0 = v(0)  \\ 
  v(y_1) &=& v_1 = \max_{y\in(0,a/2)}|v(y)| \\
  v(y_2) &=& v_2 = \min_{y\in(a/2,a)}|v(y)|. 
\end{eqnarray*}
In the case where there is not a unique point in $(0,a/2)$ at which 
$|v|$ attains its maximum $y_1$ can be chosen to be any point
at which the maximum is achieved, and similarly for $y_2$.
   
One has the obvious lower bound on the kinetic energy in terms of the 
$v_i,y_i$: 
\begin{eqnarray*}
 \int_{0}^{a} \frac{1}{2}v_{y}^2 & \ge &  \int_{0}^{y_1} \frac{1}{2}v_{y}^2 +  \int_{y_1}^{y_2} \frac{1}{2}v_{y}^2 \\
& \ge & \frac{(v_1-v_0)^2}{a}+\frac{(v_2-v_1)^2}{2a},
\end{eqnarray*}
as well as a bound on the potential energy term,
\begin{eqnarray*}
 \int_{0}^{a} Q(y)v^2 & = &  \int_{0}^{\frac{a}{2}} Q(y)v^2 +  \int_{\frac{a}{2}}^{a} Q(y)v^2 \\
 & \ge & \frac{-q_0v_1^2a}{2}+\frac{q_1v_2^2a}{2}. 
\end{eqnarray*}
The kinetic energy bound is clearly not sharp, and can be improved, 
but it suffices to prove the lemma. 
Combining these two lower bounds we find that the functional is bounded below 
by a quadratic form in three unknowns,$v_{0,1,2}$:
\begin{eqnarray*}
   \int_{-a}^{a} \frac{1}{2}v_{y}^2+Q(y)v^2 & \ge &  
   \frac{(v_1-v_0)^2}{a}+\frac{(v_2-v_1)^2}{2a} - \frac{q_0v_1^2a}{2}+\frac{q_1v_2^2a}{2}.
\end{eqnarray*}
The quadratic form is given by $v^T{\bf A}v,$ where ${\bf A}$ is defined by,
\[{\bf A} \equiv \left[\begin{array}{ccc}
\frac{1}{a}  & -\frac{1}{a} & 0 \\
-\frac{1}{a} & \frac{3}{2a}-\frac{aq_0}{2} & -\frac{1}{2a} \\
0 & -\frac{1}{2a} & \frac{1}{2a}+\frac{aq_1}{2}  \\
\end{array}\right]\]
A symmetric matrix ${\bf A}$ is positive definite if 
all of the principle minors are positive, 
and therefore $v^T{\bf A}v > 0$ when the following two inequalities are satisfied. 
\begin{eqnarray}
q_0 < \frac{1}{a^2} \label{q_0} \\
q_1-q_0-a^2q_0q_1 > 0 \label{q_1}.
\end{eqnarray}
For a fixed $a>0$ the above inequalities always have a solution 
in the positive quadrant of the $(q_0,q_1)$ plane above the hyperbola defined by $q_1-q_0-a^2q_0q_1 = 0$.For purposes of this paper we can choose 
any constants $q_0,q_1,a$ fixed (independent of $L$) such that the 
above conditions are satified. \\
 $\clubsuit$

Next, we show that we can construct a modified potential $\tilde Q$ 
such that $\tilde q(y) = \tilde Q(y)/y^2$
is smooth and the quadratic form does not decrease. This is the content of 
the next lemma.

\indent \begin{lemma} \label{SMOOTHlemma}
For any constant $\mu$ there exists a potential function $\tilde Q$ 
such that 
\begin{itemize}
\item $\tilde q = \tilde Q(y)/y^2 \in C_0^\infty$ and 
\item $\int \tilde q  \leq -\mu$
\item $\int \frac{1}{2} v_{y}^2 + \tilde Q v^2 \ge 0$  
\end{itemize}
\end{lemma}

\textbf{Proof:}
In the standard way we define $f(y)$ to be a non-decreasing $C^\infty$ function satisfying 
\[
f(y) = \left\{ \begin{array}{c} 0 ~~~ y \leq 0 \\
1 ~~~y\ge 1 \end{array} \right.
\]
Clearly we have $\lim_{y \rightarrow 0} f^{(n)}(y) = 0$ and thus 
$\lim_{y \rightarrow 0} y^{-k}f(y) =0$
Define $\tilde Q$ to be even $\tilde Q(y)=\tilde Q(-y)$ and defined for 
$y\ge0$ by
\[
\tilde Q(y) = \left\{\begin{array}{l} -q_0 f(\frac{y}{\delta}) \qquad\qquad\qquad\qquad y \in (0,\delta) \\
-q_0 \qquad\qquad\qquad\qquad\qquad y \in (\delta,\frac{a}{2} - \delta) \\
-q_0 + (q_0+q_1) f(\frac{y - \frac{a}{2} + \delta}{\delta}) ~~~ y \in (\frac{a}{2} - \delta, \frac{a}{2}) \\
q_1 \qquad\qquad\qquad\qquad\qquad\quad\!\! y \in (\frac{a}{2},a) \\
q_1 f(1 + \frac{a-y}{\delta})\qquad\qquad\qquad  y \in (a, a + \delta),\end{array}\right.
\]
where $\delta$ is small and will be chosen later. 
Clearly $\tilde q=\tilde Q/y^2$ is in $C^\infty_0$, and $\tilde q(0)=0$. 
It is easy to check that $\int \tilde q \approx -\frac{q_0}{\delta}$ 
for $\delta$ 
small, and thus the mean  $\int \tilde q $ can be made arbitrarily 
negative. 
Finally note that $\tilde Q \ge Q$ and thus  $\int \frac{1}{2} v_{y}^2 + \tilde Q v^2 \ge 0.$\\
 $\clubsuit$

We are now in a position to prove our main theorem, in which 
we construct our potential $\phi_x.$ The Cauchy-Schwartz inequality implies
that $-\int u_x^2 \ge -\frac{1}{2}\int u_{xx}^2 - \frac{1}{2} \int u^2$, so 
it clearly suffices to show that 
$$\int \frac{1}{4} u_{xx}^2 + (\phi_x - \frac{3}{4}) u^2 \ge 0 $$   
We write the mean zero function $\phi_x$ in the form  
\[
\phi_x = q(x) - <\!\!q\!\!>
\]
where $<\!\cdot\!>$ denotes the mean value on $[-L,L]$:
\[
<\!\! q \!\!> = \frac{1}{2L} \int_{-L}^L q(x) dx.
\]
If we can choose $q$ such that $<\!\!q\!\!>\le -\frac{3}{4}$ and 
$$
\int \frac{1}{4} u_{xx}^2 + q(x) u^2 \ge 0
\label{QFP}
$$
then we are done. 

Based on the scaling arguments of the previous section 
we define $q(x)$ to be 
\[
q(x) = L^\frac{4}{3} \tilde q(x L^\frac{1}{3})  ~~~~x \in [-L,L]
\]
where $\tilde q$ is compactly supported. 
Obviously this can be extended to a periodic function 
in the standard way. After the rescaling $y = L^\frac{1}{3} x$ the quadratic 
form in equation (\ref{QFP}) becomes 
\[
L \int_{-L^{\frac{4}{3}}}^{-L^{\frac{4}{3}}}
\frac{1}{4} u_{yy}^2 + \tilde q(y) u^2 dy.  
\]
From lemma \ref{hardylemma} we have the inequality
\[
\int \frac{1}{4} u_{yy}^2 + \tilde q(y) u^2 dy \ge \int \frac{1}{2} v_{y}^2 + 
y^2 \tilde q(y) v^2 dy = \int \frac{1}{2} v_{y}^2 + \tilde Q(y) v^2 dy 
\]
From the results of lemma \ref{SMOOTHlemma} we can choose $\tilde Q$ in such
a way that the the above is positive 
for all $u \in C_0^\infty$ with $u(0)=0$. Since this is a dense domain 
the operator is positive under Dirichlet boundary conditions. 

\begin{remark} 
The fact that $\tilde q$ can be chosen to have arbitrarily 
negative mean and still generate a positive operator $K$ implies that 
$\phi_x$ can be chosen such that $<\!u,Ku\!> \ge \delta |\!|u|\!|^2$ for 
any constant $\delta$ independent of $L$ with a bound 
$|\!|\phi|\!|_{H^2}\le C L^\frac{3}{2},$ where $C$ depends on $\delta$ 
but is independent of $L$. Thus the above argument gives an
$L_2$ of the destabilized KS equation
\[
u_t = -u_{xxxx} - u_{xx} + \gamma u + u u_x 
\] 
which scales like $L^\frac{3}{2}$ for any fixed $\gamma$. 
\end{remark}

We have constructed a potential function $\phi_x$ such that the 
operator ${K}$ is bounded below on the set of functions 
satisfying a Dirichlet boundary condition. All that remains to be 
checked is that the $H^2$ norm of the potential scales correctly with $L$. 
This is the content of the next lemma.

\indent \begin{lemma} \label{biglemma}
The potential $\phi$ satisfies $|\!|\phi|\!|_{H^2} \le C L^\frac{3}{2}$.
\end{lemma}
{\bf Proof:} From the definition 
\[
\phi_x = <\!\!q\!\!> - q(x) = <\!\!q\!\!> - L^\frac{4}{3}\tilde q(xL^\frac{1}{3}) 
\]
it is clear on rescaling that $|\!|\phi_x|\!|_2^2 = O(L^\frac{7}{3})$ and that 
 $|\!|\phi_{xx}|\!|_2^2 = O(L^3)$. Thus we only need to estimate 
 $|\!|\phi|\!|_2^2$. 
Since $\phi_x$ is defined to be 
\begin{eqnarray*}
\phi_x(x) = <\!\!q\!\!> - q(x) 
\end{eqnarray*}
we have 
\begin{eqnarray*}
\phi(x) = \int_0^x \phi_s ds = \int_0^x q(s) - <\!\!q\!\!>ds,
\end{eqnarray*}
and after the substitution $y = sL^{1/3}$ this becomes
\[
\phi(x) = L \int_0^{xL^{1/3}} \tilde q dy - <\!\!q\!\!>x.
\]
We have the obvious estimate 
\[
|\phi(x)| \le L\int_0^{xL^{1/3}} |\tilde q(y)| dy + <\!\!q\!\!> L
\]
Since $\tilde q$ is bounded (independently of $L$) and supported on 
$[-a,a]$ we have 
\[
|\phi(x)| \le c a L + <\!\!q\!\!> L = O(L)
\]
The $L_2$ bound now follows since
\[ |\!|\phi|\!|_2^2 \le (2L)|\!|\phi|\!|_{\infty}^2 = O(L^3) 
\]
$\clubsuit$

\subsection{Extension to Arbitrary Initial Data}

The theorem proved above, together with the lemmas proved in the first section, 
show that for odd initial data the (destabilized) Kuramoto-Sivashinsky equation 
remains bounded in $L_2$ for all time.  
These results can be extended to arbitrary mean-zero initial data 
in the manner first done by Collet, Eckmann, Epstein and Stubbe\cite{CEES}
or by Goodman\cite{Good}. This paper was written in such a way as to be 
compatible with the results of Collet et. al, where one allows the 
potential $\phi_x$ to translate via a gradient-flow type dynamics. 
In particular Theorem 1 in this paper is 
stated in such as way as to be compatible with Lemma 5.1 in \cite{CEES}. 
From this the results of Proposition 4.3 in \cite{CEES} follow, and the 
$L_2$ boundedness result extends to arbitrary mean-zero initial data. 
One could equally well employ the related idea of Goodman, and look 
at the rate of change of the distance of $u$ from the set of all 
translates of $\phi$.

\section{CONCLUSIONS}

In this paper we have constructed a function $\phi$ such that the ball 
$B(\phi,c L^\frac{3}{2})$ is a global attracting set for the (destabilized) 
Kuramoto-Sivashinsky equation 
$$
u_t = -u_{xxxx} - u_{xx} + \gamma u + u u_x. 
$$
This result has the best possible scaling with $L$, the size of the domain, since 
the above equation has stationary solutions with $L_2$ norm which scales like 
$L^\frac{3}{2}$. While the result of Giacomelli and Otto is slightly 
stronger, the 
Lyapunov function argument outlined here is still interesting for a number of 
reasons. First is that it gives the optimal scaling for $L_2$ boundedness 
of the 
destabilized KS equation, a result that seems to be new. Secondly, this 
calculation 
shows that the critical exponents for an argument of this type can be 
achieved, and 
is interesting as a demonstration of the limits of this kind of Lyapunov 
function 
argument. Finally, it should be noted that there is an elegant
 result of Molinet\cite{Mol}, improving on earlier work of Sell and 
Taboada\cite{ST}, 
which gives $L_2$ boundedness of the Kuramoto-Sivashinsky equation in two 
spatial dimensions 
$$
\vec u_t = -\Delta^2 \vec u - \Delta \vec u + \vec \nabla 
(\vec u \cdot \vec u) ~~~~~~~~~~ \vec \nabla^\perp \cdot \vec u = 0  
$$
for sufficiently thin rectangular domains. The results of Molinet require 
the construction of a 
Lyapunov function for the problem one spatial dimension. In the original paper 
Molinet uses the 
Lyapunov function constructed by Collet et. al. to show boundedness in 
$L_2[(0,L_x) \times (0,L_y)]$ 
assuming that the width in the second spatial direction satisfies (for $L_x >\!> 1,L_y <\!< 1$) 
\[
L_y \leq C L_x^{-\frac{67}{35}}.
\]
Assuming the above condition holds, together with a bound on the 
$L_2$ of the initial data,  Molinet establishes a long-time bound on the $L_2$ norm of the form
\[
\limsup _{t\rightarrow \infty}||\vec u||_2 \le C L_x^\frac{8}{5} L_y^\frac{1}{2}.
\]
In fact Molinet shows a great deal more, including decay in time of the 
second component of $\vec u$ and explicit estimates of the relevant constants. 
Using the Lyapunov function constructed here and applying Molinet's results verbatim 
gives $L_2$ boundedness assuming the aspect ratio satisfies 
\[
L_y \leq C L_x^{-\frac{13}{7}}.
\]
This, in turn, leads to a bound on the $L_2$ norm of the form
\[
\limsup_{t\rightarrow\infty}||\vec u||_2 \le C L_x^\frac{3}{2} L_y^\frac{1}{2}.
\]

{\bf Acknowledgements:} The authors would like to thank Dirk Hundertmark, Ralf Wittenberg and Jonathan Goodman for many useful 
conversations. In particular we would like to thank Ralf Wittenberg for 
pointing out to us the significance of the destabilized Kuramoto-Sivashinsky 
equation, and the connection to the background flow method. 
This research was supported by NSF grants DMS-0203938 and DMS-0354373. 
JCB would like to thank 
the Stanford University mathematics department for hospitality during the writing of part of this paper.

\end{document}